\input amstex
\documentstyle{amstexsi}
\firstpageno{1}                                                   %
\shortauthor{Christopher J. Hillar and Charles R. Johnson}        %
\shorttitle{Eigenvalues of Generalized Words}                     %
      %
 %
 %
         %
           %
\topmatter
\title
Positive Eigenvalues of Generalized Words in Two Hermitian
Positive Definite Matrices\footnote[\boldkey*]{This research was
conducted, in part, during the summer of 1999 at the College of
William and Mary's Research Experiences for Undergraduates program
and was supported by NSF REU grant DMS-96-19577.}
\endtitle
\author
Christopher J. Hillar\footnote[\dag]{Department of Mathematics,
University of California, Berkeley, CA 94720. \\
(chillar\@math.berkeley.edu). The work of the first author is
supported under a National Science Foundation Graduate Research
Fellowship}\ and Charles R. Johnson\footnote[\ddag]{Department of
Mathematics, College of William and Mary, Williamsburg, VA
23187-8795. (crjohnso\@math.wm.edu)}
\endauthor
\abstract We define a word in two positive definite (complex
Hermitian) matrices $A$ and $B$ as a finite product of real powers
of $A$ and $B$. The question of which words have only positive
eigenvalues is addressed. This question was raised some time ago
in connection with a long-standing problem in theoretical physics,
and it was previously approached by the authors for words in two
real positive definite matrices with positive integral exponents.
A large class of words that do guarantee positive eigenvalues is
identified, and considerable evidence is given for the conjecture
that no other words do.
\endabstract
\endtopmatter
\document
\subhead 1. Introduction\endsubhead A {\it generalized word} ({\it
g-word}, for short) $W = W(A,B)$ in two letters $A$ and $B$ is an
expression of the form $$ W = A^{p_1 } B^{q_1 } A^{p_2 } B^{q_2 }
\cdots A^{p_k } B^{q_k } A^{p_{k + 1} } $$ in which the exponents
$p_{i}$ and $q_{i}$ are real numbers such that $p_{i},q_{i} \neq
0$, $i = 1,\ldots,k$, and $p_{k+1}$ is an arbitrary real number.
We call $k$ the {\it class number} of $W$. The {\it reversal} of
the g-word $W$ is $W^{*} = A^{p_{k + 1} } B^{q_k } A^{p_k } \cdots
B^{q_2 } A^{p_2 } B^{q_1 } A^{p_1 }$ and a g-word is {\it
symmetric} if it is identical to its reversal (in other contexts,
the name ``palindromic" is also used).

We are interested in the matrices that result when the two letters
are (independent) positive definite (complex Hermitian) $n$-by-$n$
matrices (PD, for short).  For convenience, the letters $A,B$ will
also represent the substituted PD matrices (the context will make
the distinction clear).  To make sure that $W$ is well-defined
after substitution, we take primary PD powers (see [4, p.433] and
[4, p.413]). I.e., given $p \in \Bbb R$, a unitary matrix $U$, and
a positive diagonal matrix $D$, we have ${(UDU^{*})}^{p} =
UD^{p}U^{*}$.

We are primarily interested in those g-words for which $W$ has
only positive (real) eigenvalues, no matter what the positive
definite matrices $A$ and $B$ are (for any positive integer $n$).
We call such a g-word {\it good} and all other g-words {\it bad}.
Interest in this problem stems from a question in quantum physics
[1,7], as discussed in [5] for the case of (ordinary) words
($p_{i},q_{i}$ positive integers) and real symmetric positive
definite matrices. For 2-by-2 matrices, the situation is better
understood.  For example, if all of the $p_{i}$ or $q_{i}$ are of
the same sign, then it is known [2] that any word in two 2-by-2 PD
matrices necessarily has positive eigenvalues. Since we are
interested in those words that are bad, our search should begin
with $n = 3$.

We call a g-word {\it nearly symmetric} if it is either symmetric
or a product (juxtaposition) of two symmetric words.  It is an
elementary exercise that good g-words (and, therefore, bad also)
are unchanged by each of the following:

\

i) reversal;
\

ii) interchange of the letters $A$,$B$;
\

iii) cyclic permutation, e.g., $$A^{p_1 } B^{q_1 } A^{p_2 } B^{q_2
} \cdots A^{p_k } B^{q_k }  \longrightarrow B^{q_1 } A^{p_2 }
B^{q_2 } \cdots A^{p_k } B^{q_k } A^{p_1 }.$$

iv) multiplication of all the $p_{i}$'s ($q_{i}$'s) by a fixed
nonzero scalar.

\

It is also an elementary combinatorial exercise that each of
(i)-(iv) preserves the nearly symmetric g-words.

Recall that two $n$-by-$n$ matrices $X$ and $Y$ are said to be
{\it congruent} if there is an invertible $n$-by-$n$ matrix $Z$
such that $Y = Z^{*}XZ$ and that congruence on Hermitian matrices
preserves inertia (the ordered triple consisting of the number of
positive, negative, and zero eigenvalues) and, thus, positive
definiteness [3, p.223].  A symmetric word of class $k$ in two
positive definite matrices is congruent to one of class $k-1$,
iteration of which implies congruence to the ``center," class 0,
positive definite matrix.  We conclude that

\

\proclaim{Lemma 1.1}\rm  A symmetric g-word in two positive
definite matrices is positive definite; thus, every symmetric
g-word is good.
\endproclaim

\

It is also known [3, p.465] that a product of two positive
definite matrices has only positive eigenvalues.  In view of Lemma
1.1, it then follows that

\

\proclaim{Theorem 1.2}\rm  Each nearly symmetric g-word is good.
\endproclaim

\

We conjecture the converse to Theorem 1.2.

\

\proclaim{Conjecture 1.3}\rm  A g-word is good if and only if it
is nearly symmetric.
\endproclaim

\

We note that near symmetry is easily verified algorithmically.
Using iii), each g-word $W$ may be taken to be in the form in
which $p_{k+1} = 0$, i.e. the g-word may be taken to begin with a
power of one letter and end with a power of the other.  We refer
to this as {\it standard form}.  Then, near symmetry of $A^{p_1 }
B^{q_1 } A^{p_2 } B^{q_2 }  \cdots A^{p_k } B^{q_k }$ may be
determined by inspection of each pair of words of the form $A^{p_1
} B^{q_1 } A^{p_2 } B^{q_2 }  \cdots A^{p_i }$, $B^{q_i } A^{p_{i
+ 1} } B^{q_{i + 1} }  \cdots A^{p_k } B^{q_k }$, $i =
1,\ldots,k$.

The number of class $k$ words that are nearly symmetric is small 
compared to the total number of all such words.  Therefore, 
if Conjecture 1.3 were to be true, it is necessary that the
``density" of good words is also 0.  To make this precise,
we view the space of good words (of class $k$) as a subset of
$\Bbb R^{2k}$ parameterized by $\{(p_{1},q_{1},...,p_{k},q_{k}) \
| \ W$ is good$\}$.  Our main theorem in this direction is 
then given by the following.

\

\proclaim{Theorem 1.3}\rm  The set of good words
of class $k$ has measure zero in $\Bbb R^{2k}$.
\endproclaim

\

Each g-word of class 0 or 1 is nearly symmetric, and a standard
class 2 g-word ($A^{p_1 } B^{q_1 } A^{p_2 } B^{q_2 }$) is nearly
symmetric if and only if $p_{1} = p_{2}$ or $q_{1} = q_{2}$.  We
first show that Conjecture 1.3 is correct for class 2 words, in
part because the technique generalizes in certain ways.  We then
describe necessary conditions for a word to be good (Theorems 3.1
and 4.1), thereby showing the rarity of such words. The idea is to
produce PD matrices $A$ and $B$ for which the trace of a given
word has nonzero imaginary part. Thus, not all eigenvalues could
be positive.

\

\subhead 2. Class 2 words\endsubhead To prepare for the class 2
discussion and the general results that follow, some preliminaries
are necessary.

\

\proclaim{Definition 2.1}\rm  A {\it generalized polynomial} is an
expression of the form

$$G(x_{1},x_{2},\ldots,x_{m})  =  \sum\limits_{i = 1}^q {c_i
\prod\limits_{j = 1}^m {x_j ^{p_{i,j} } } }$$ in which $c_{i} \in
\Bbb C, p_{i,j} \in \Bbb R$, and the $x_{j}$ are the variables.
\endproclaim

\

A generalized polynomial ({\it g-poly}, for short) is said to be
{\it reduced} if for each $i \neq t$, there is a $j$ such that
$p_{i,j} \neq p_{t,j}$. For example, $2x_{1}^{0}x_{2}^{0} +
x_{1}^{-2}x_{2} + 3x_{1}x_{2}^{0}$ is a reduced g-polynomial, but
$x_{1}^{2}x_{2} - 3x_{1}x_{2} + 5x_{1}^{2}x_{2}$ is not.

The proof in the class 2 case uses expressions that are
generalized polynomials.  They are also important in the study of
larger class numbers.  Let $\Bbb R_{+}^{m} =
\{(x_{1},x_{2},...,x_{m}) \ | \ x_{j} > 0, \ j = 1,\ldots,m\}$ and
let $$R_{G} = \{ (x_{1},x_{2},...,x_{m}) \in \Bbb R_{+}^{m} \ | \
G(x_{1},x_{2},...,x_{m}) = 0\},$$ the set of positive zeroes of
$G$.  We first write down a preliminary lemma for the one variable
case.

\

\proclaim{Lemma 2.2}\rm  If $G$ is a reduced 1-variable ($m=1$)
g-poly and if 0 is an accumulation point for $R_{G}$, then $c_{i}
= 0$ for all $i$.
\endproclaim

\

\demo{Proof} Assume $G$ is reduced with one variable $x_{1}$ and
for some $i \in \{1,\ldots,q\}$, $c_{i} \neq 0$.  Also, suppose 0
is an accumulation point for $R_{G}$.  Multiplying $G$ by a large
enough power of $x_{1}$, we can assume that $G$ has only positive
exponents -- since this does not change $R_{G}$. Let $M = \mathop
{\min }\limits_i \{p_{i,1}\}$, and notice that $$ G = x_1 ^M
\left( {\sum\limits_{i = 1}^q {c_i \,x_1 ^{p_{i,1}  - M} } }
\right) = x_1 ^M G_2.$$

Since $G$ is reduced, so is $G_{2}$, and $G_{2}$ has a non-zero
constant term. Also, since $G_{2}$ has the same set of positive
zeroes as $G$, it follows that $R_{G_2 } = R_G$.  Now, 0 is an
accumulation point for $R_{G}$, so there exists a sequence of
$\lambda_{r} \in R_{G}$ ($r = 1,2,\ldots$) such that $\mathop
{\lim }\limits_{r \to \infty } \lambda_r = 0$. Since,
$G_{2}(\lambda_{r}) = 0$ for all $r$, we have that $0 = \mathop
{\lim }\limits_{r \to \infty } G_2 (\lambda _r )$. But $G_{2}$ is
continuous, and so $G_{2}(0) = 0$. This contradicts the fact that
$G_{2}$ has a non-zero constant term, proving the lemma.
\qquad\qed
\enddemo

\

The lemma above allows us to prove the main observation we need
about g-polys.  We remark that this result is similar to one in
[6, p.176] concerning (normal) polynomials over a field.

\

\proclaim{Lemma 2.3}\rm  Let $S_{1},\ldots,S_{m}$ be infinite
subsets of $\Bbb R_{+}$ with 0 being an accumulation point for
each. If $G$ is a reduced g-poly and if $G(a_{1},\ldots,a_{m}) =
0$ for all $a_{j} \in S_{j} \ (j=1,\ldots,m)$, then $c_{i} = 0$
for all $i$.
\endproclaim

\

\demo{Proof} For $m = 1$ and $q$ arbitrary, we have Lemma 2.2
above. Therefore, we consider induction on $m$.  Assume $G$ is a
reduced g-poly and for some $i \in \{1,\ldots,q\}$, $c_{i} \neq
0$. Multiplying $G$ by a large enough power of $\prod\nolimits_{j
= 1}^m {x_j }$, we can assume that $G$ has only positive
exponents, as this does not change $R_{G}$. Let $M = \mathop {\min
}\limits_i \{p_{i,m}\}$, and examine $$ G = x_m ^M \left(
{\sum\limits_{i = 1}^q {c_i \,x_m ^{p_{i,m} - M} \prod\limits_{j =
1}^{m - 1} {x_j ^{p_{i,j} } } } } \right) = x_m ^M G_2.$$

Since $G$ is reduced, so is $G_{2}$, and $G_{2}$ has non-zero
terms that do not contain the variable $x_{m}$ (i.e. those terms
that only contain $x_{m}^{0}$). Specifically, let
$F(x_{1},\ldots,x_{m-1})$ be the reduced g-poly in $G_{2}$ that
does not contain the variable $x_{m}$. Since $G_{2}$ has the same
set of positive zeroes as $G$, it follows that $R_{G_2} = R_G$.
Now, fix $(a_{1},\ldots,a_{m-1})$ $\in$ $S_{1}\times \cdots \times
S_{m-1}$ and examine the g-poly in 1 variable,
$G_{2}(a_{1},\ldots,a_{m-1},x_{m})$.  By Lemma 2.2, we must have
that the constant term of this g-poly is zero.  Whence, for all
$(a_{1},\ldots,a_{m-1}) \in S_{1}\times \cdots \times S_{m-1}$, it
follows that $F(a_{1},\ldots,a_{m-1}) = 0$.  By induction, each
term in $F$ must have zero coefficient, and this contradiction
finishes the proof.  \qquad\qed
\enddemo

\

For the rest of the discussion, we will assume the following
parameterization of $A$ and $B$.  Let $S = UU^{T}$ for some
unitary matrix $U$, so that $S$ is both unitary and symmetric.
Then, we will assume $$S = \left( {\matrix
   {s_{1,1} } & {s_{2,1} } & {s_{3,1} }  \\
   {s_{2,1} } & {s_{2,2} } & {s_{3,2} }  \\
   {s_{3,1} } & {s_{3,2} } & {s_{3,3} }  \\
 \endmatrix } \right), \ A = \left( {\matrix
   1 & 0 & 0  \\
   0 & x_{1} & 0  \\
   0 & 0 & y_{1}  \\
 \endmatrix } \right),  \ E = \left( {\matrix
   1 & 0 & 0  \\
   0 & x_{2} & 0  \\
   0 & 0 & y_{2}  \\
 \endmatrix } \right), \ B = SE\overline S \tag 2.1$$ for  $x_{1},x_{2},y_{1},y_{2} > 0$.  We can view the trace of the g-word,
$W$, under the assumption of (2.1) as $$ \text{Tr}\left[ {A^{p_1 }
SE^{q_1 } \overline S A^{p_2 } SE^{q_2 } \overline S \cdots A^{p_k
} SE^{q_k } \overline S } \right].$$

We are now ready to prove

\

\proclaim{Theorem 2.4}\rm  A class 2 g-word is good if and only if
it is nearly symmetric.
\endproclaim

\

\demo{Proof} We may suppose that $W = A^{p_1 } B^{q_1 } A^{p_2 }
B^{q_2 }$ is a class 2 g-word in standard form.  If $W$ is nearly
symmetric, Theorem 1.2 shows that $W$ is good.  If $W$ is not
nearly symmetric, then we find a unitary matrix $U$ and positive
diagonal matrices $A$ and $E$ such that $B = SE\overline S$ gives
$W$ a non-real trace.  By the above remark, the trace of $W$ is
that same as the trace of $A^{p_1 } SE^{q_1 } \overline S A^{p_2 }
SE^{q_2 } \overline S$ .  Set $x = x_{1} = x_{2}$, $y = y_{1} =
y_{2}$ and consider $A =$ diag$(1,x,y) = E$.  Also, take $U$ to be
the unitary matrix $$U = \frac{1} {4}\left( {\matrix
   2 & { - 1 - i} & {3 - i}  \\
   { - 1 + i} & 3 & {1 - 2i}  \\
   {3 + i} & {1 + 2i} & { - 1}  \\
 \endmatrix } \right).$$

With these assumptions, a straightforward computation reveals that
the imaginary part of Tr[$A^{p_1 } SA^{q_1 } \overline S A^{p_2 }
SA^{q_2 } \overline S$] is a g-poly in $x$ and $y$ given by
$\frac{3} {{128}}$ times $$\left[ {x^{p_2 } \left( {y^{p_1 }  - 1}
\right) + x^{p_1 } \left( {1 - y^{p_2 } } \right) + y^{p_2 }  -
y^{p_1 } } \right] \cdot \left[ {x^{q_2 } \left( {y^{q_1 }  - 1}
\right) + x^{q_1 } \left( {1 - y^{q_2 } } \right) + y^{q_2 }  -
y^{q_1 } } \right].$$

Now, fix $y > 1$ and suppose the above expression is 0 for all $x
> 0$ (which is implied if $W$ is good).  Then, one of the factors
has an accumulation point of 0 in its set of positive zeroes.
Therefore, from Lemma 2.2, one of these factors is not reduced.
In this case, we must have either $p_{2} = p_{1}$ or $q_{2} =
q_{1}$. This completes the proof of the theorem.    \qquad\qed
\enddemo

\

\subhead 3. Good words are rare\endsubhead We next show that the
``density" of good words is 0, partly extending Theorem 2.4 in the
direction of a proof of our conjecture for higher class words. We
will prove that for good $W$ there are certain non-trivial
algebraic (actually linear) relations linking the $p_{i},q_{i}$.
If we view the space of good words (of class $k$) as a subset of
$\Bbb R^{2k}$ parameterized by $\{(p_{1},q_{1},...,p_{k},q_{k}) \
| \ W$ is good$\}$, then it will follow that the set of good words
is a set of measure zero in $\Bbb R^{2k}$. Consider a subset $P$
of $\{p_{1},\ldots,p_{k}\}$ (as variables) and similarly let $Q$
be a subset of $\{q_{1},\ldots,q_{k}\}$ (also as variables). Then,
the full statement is given by

\

\proclaim{Theorem 3.1}\rm  If $W = A^{p_1 } B^{q_1 } A^{p_2 }
B^{q_2 }  \cdots A^{p_k } B^{q_k }$ is good, then there exists a
subset, $P \subseteq \{p_{1},\ldots,p_{k}\}$, and a subset, $Q
\subseteq \{q_{1},\ldots,q_{k}\}$, such that $$p_{1} =
\sum\limits_{p \in P }{p } \ \ \ \text{and} \ \ \ q_{1} =
\sum\limits_{q \in Q }{q},$$ in which this pair of relations is
not the trivial one, $$p_{1} = p_{1} \ \ \ \text{and} \ \ \ q_{1}
= q_{1}.$$
\endproclaim

\

\proclaim{Corollary 3.2}\rm If the $p_{i}$ are linearly
independent over $\{-1,0,1\}$ and if the $q_{i}$ are also linearly
independent over $\{-1,0,1\}$, then $W$ is bad.
\endproclaim

\

As another immediate Corollary, we obtain Theorem 1.3 mentioned 
in the introduction.

We should remark that Theorem 3.1 gives another proof of Theorem
2.4.  For in the case of class 2 words, we must have (since the
$p_{i},q_{i} \neq 0$) either $p_{1} = p_{2}$ or $q_{1} = q_{2}$
for $W$ to be good.  Before proving this theorem we need a few
technical remarks. Assuming the parameterization as in (2.1), the
trace of $W$ can be viewed as a g-poly in $x_{1}$, $x_{2}$,
$y_{1}$, and $y_{2}$ with exponents involving the $p_{i}$ and
$q_{i}$ : $$ \text{Tr}\left[ {A^{p_1 } SE^{q_1 } \overline S
A^{p_2 } SE^{q_2 } \overline S \cdots A^{p_k } SE^{q_k } \overline
S } \right].$$  In fact, it is not hard to see that this trace
will be a constant plus a (formal) sum of elements of the form, $$
c\,x_{1}^{\sum\limits_{p \in P_{1} }{p}} y_{1}^{\sum\limits_{p \in
P_{2} }{p}}x_{2}^{\sum\limits_{q \in Q_{1} }{q}}
y_{2}^{\sum\limits_{q \in Q_{2} }{q}}, \tag3.1$$ in which $c \in
\Bbb C$, $P_{1},P_{2} \subseteq \{p_{1},\ldots,p_{k}\}$, and
$Q_{1},Q_{2} \subseteq \{q_{1},\ldots,q_{k}\}$. Of course, here an
empty sum is defined to be zero. We call an expression as above a
{\it term} in the trace expansion of $W$. We will also say that
$p_{i}$ ($q_{i}$) is {\it contained} in a power of $x_{1}$,
$y_{1}$ ($x_{2}$, $y_{2}$) if $p_{i} \in P_{1}, p_{i} \in P_{2}$
($q_{i} \in Q_{1}, q_{i} \in Q_{2}$).  As an example, consider
Tr[$A^{p_{1}}SE^{q_{1}}\overline{S}$], corresponding to the g-word
$W = A^{p_{1}}B^{q_{1}}$.  This trace is $$s_{11} \overline
{s_{11} } + s_{21} \overline {s_{21} } x_2 ^{q_1 }  + s_{31}
\overline {s_{31} } y_2 ^{q_1 }  + s_{21} \overline {s_{21} } x_1
^{p_1 }  + s_{22} \overline {s_{22} } x_1 ^{p_1 } x_2 ^{q_1 }  +
s_{32} \overline {s_{32} } x_1 ^{p_1 } y_2 ^{q_1 }$$ $$+ s_{31}
\overline {s_{31} } y_1 ^{p_1 }  + s_{32} \overline {s_{32} } x_2
^{q_1 } y_1 ^{p_1 } + s_{33} \overline {s_{33} } y_2 ^{q_1 } y_1
^{p_1 }.$$ Since these expressions are quite complicated (even for
$k = 1$), we need a method of isolating the coefficients of
individual terms as functions of the entries of $S$.  The
following lemma describes a way to determine these quantities.

\

\proclaim{Lemma 3.3}\rm Assume we have a representation (2.1).
Consider an arbitrary term as in (3.1), and let $P = P_{1} \cup
P_{2}$ $= \{ p_{i_1 },\ldots ,p_{i_u }\}$ and $Q = Q_{1} \cup
Q_{2}$ $= \{ q_{j_1 } ,\ldots ,q_{j_v } \}$.  So each element of
$P$ and $Q$ is contained in a power of $x_{1}$, $x_{2}$, $y_{1}$,
or $y_{2}$ in the given term of Tr$[W] =$ Tr$[A^{p_1 } SE^{q_1 }
\overline S A^{p_2 } SE^{q_2 } \overline S \cdots A^{p_k } SE^{q_k
} \overline S]$. Then, this term will appear with the same
coefficient in the trace of:

$$ A^{p_{i_1 } } SF\overline S FSF\overline S  \cdots SE^{q_{j_1 }
} \overline S FS \cdots SE^{q_{j_v } } \overline S  \cdots
SF\overline S $$ where $F = \left( {\matrix
   1 & 0 & 0  \\
   0 & 0 & 0  \\
   0 & 0 & 0  \\
 \endmatrix } \right)
$ replaces each instance of $A^{p_i }$ ($E^{q_j }$) in which
$p_{i} \notin P$ ($q_{j} \notin Q$).
\endproclaim

\

\

\demo{Proof} Let $X_{1,1} = x_1 ^{p_1 },\ldots,X_{1,k} = x_1 ^{p_k
}$; $X_{2,1} = x_2 ^{q_1 },\ldots,X_{2,k}  = x_2 ^{q_k }$, and
similarly, let $Y_{1,1} = y_1 ^{p_1 },\ldots,Y_{1,k} = y_1 ^{p_k
}$; $Y_{2,1} = y_2 ^{q_1 },\ldots,Y_{2,k} = y_2 ^{q_k }$. Now, set

$$A_{i} = \left( {\matrix
   1 & 0 & 0  \\
   0 & {X_{1,i} } & 0  \\
   0 & 0 & {Y_{1,i} }  \\
 \endmatrix } \right), \ \ E_{i} = \left( {\matrix
   1 & 0 & 0  \\
   0 & {X_{2,i} } & 0  \\
   0 & 0 & {Y_{2,i} }  \\
 \endmatrix } \right)$$ for $i = 1,\ldots,k$ and examine the trace of $$H = A_1 SE_1 \overline S A_2 SE_2 \overline S  \cdots A_k SE_k
\overline S.$$ This trace is a (normal) polynomial in $X_{1,i}$,
$X_{2,i}$, $Y_{1,i}$, and $Y_{2,i}$. Next, consider the trace of
$$J = A_{i_1 } SF\overline S FSF\overline S  \cdots SE_{j_1 }
\overline S FS \cdots SE_{j_v } \overline S  \cdots SF\overline S
,$$ where $F$ as above replaces in $H$ each instance of $A_{i}$
($E_{j}$) in which $p_{i} \notin P$ ($q_{j} \notin Q$).  The trace
of $J$ is also a polynomial in
$X_{1,i}$,$X_{2,i}$,$Y_{1,i}$,$Y_{2,i}$.  But now, notice that if,
for example, $p_{i} \notin P$, then every time a factor (when
computing this trace) of $X_{1,i}$ would have appeared in Tr$[H]$,
it is replaced by a 0 -- similarly for $X_{2,i}$,$Y_{1,i}
$,$Y_{2,i}$. Thus, it is clear that Tr$[J]$ has the desired term
with the same coefficient as in Tr$[H]$. \qquad\qed
\enddemo

\

Our first application of this result is the following.

\

\proclaim{Lemma 3.4}\rm Assuming that we have representation
(2.1), the g-poly expression, Tr$[A^{p_1 } SE^{q_1 } \overline S
A^{p_2 } SE^{q_2 } \overline S  \cdots A^{p_k } SE^{q_k }
\overline S]$, has a real constant term.
\endproclaim

\

\demo{Proof} Using Lemma 3.3, the constant term of the above trace
is just the trace of: $ \left( {FSF\overline S } \right)\left(
{FSF\overline S } \right) \cdots \left( {FSF\overline S } \right)
= \left( {FSF\overline S } \right)^k$. But $$ \left( {FSF\overline
S } \right)^k = \left( {\matrix
   {s_{1,1} \overline {s_{1,1} } } & {s_{1,1} \overline {s_{2,1} } } & {s_{1,1}
   \overline {s_{3,1} } }  \\
   0 & 0 & 0  \\
   0 & 0 & 0  \\
 \endmatrix } \right)^k
$$

$$ = \left( {\matrix
   {\left( {s_{1,1} \overline {s_{1,1} } } \right)^k } &
   {s_{1,1} \overline {s_{2,1} } \left( {s_{1,1} \overline {s_{1,1} } }
   \right)^{k - 1} } & {s_{1,1} \overline {s_{3,1} } \left( {s_{1,1} \overline {s_{1,1} } }
   \right)^{k - 1} }  \\
   0 & 0 & 0  \\
   0 & 0 & 0  \\
 \endmatrix } \right).
$$  Hence the trace of $A^{p_1 } SE^{q_1 } \overline S \cdots
A^{p_k } SE^{q_k } \overline S$ has a real constant term, $\left(
{s_{1,1} \overline {s_{1,1} } } \right)^k$. \qquad\qed
\enddemo

\

Call two exponents $p_{i}$ and $q_{j}$ {\it adjacent} if either $i
= j$, $i - 1 = j$, or $(i,j) = (1,k)$.  For instance, $p_{2}$ and
$q_{1}$ are adjacent.  The following computation shows that we can
always find terms with non-real coefficients.

\

\proclaim{Lemma 3.5}\rm Assuming (2.1) as before, the coefficient
of any term, $x_{1}^{p_{i_1}} y_{2}^{q_{j_1}}$, in the trace of a
class $k$ word ($k > 1$) in which $p_{i_1}$ and $q_{j_1}$ are
adjacent is given by:

$$ \overline {s_{2,1} } s_{1,1} \overline {s_{3,1} } s_{3,2}
\left( {s_{1,1} \overline {s_{1,1} } } \right)^{k - 2} \ \
\text{or} \ \ \ s_{2,1} \overline {s_{1,1} } s_{3,1} \overline
{s_{3,2} } \left( {s_{1,1} \overline {s_{1,1} } } \right)^{k -
2}.$$
\endproclaim

\

\demo{Proof} First notice that we can assume we are dealing with
the term $x_{1}^{p_1 } y_{2}^{q_1 }$ by (possibly) a reversal and
(possibly) cycling.  (A cycling is a similarity transformation,
and a reversal corresponds to a conjugate transposition, which
will only change the conjugacy of the coefficient). Therefore, we
need only compute (by Lemma 3.3) $$ A^{p_1 } SE^{q_1 } \overline S
FSF\overline S \cdots FSF\overline S.$$ This quantity can be
represented compactly as: $$ A^{p_1 } SE^{q_1 } \overline S \left(
{FSF\overline S } \right)^{k - 1}.$$ A straightforward calculation
gives us that the trace of the product of the two matrices $A^{p_1
} SE^{q_1 } \overline S$ and $\left( {FSF\overline S } \right)^{k
- 1}$ produces a coefficient of $$ \overline {s_{2,1} } s_{1,1}
\overline {s_{3,1} } s_{3,2} \left( {s_{1,1} \overline {s_{1,1} }
} \right)^{k - 2} $$ for the $x_{1}^{p_1 } y_{2}^{q_1 }$ term as
stated (and this term appears only once).  Notice that this
coefficient can be made non-real (for instance, using $U$ as in
the proof of Theorem 2.4). \qquad\qed
\enddemo

\

We are now ready to prove Theorem 3.1.

\

\demo{Proof of Theorem 3.1} Assuming the parameterization (2.1),
notice that the imaginary part of Tr$[W]$ is a g-poly in $x_{1}$,
$x_{2}$, $y_{1}$, and $y_{2}$ with exponents involving the $p_{i}$
and $q_{j}$.  This g-poly has a term $x_{1}^{p_1 } y_{2}^{q_1 }$
with a nonzero coefficient by Lemma 3.5, and it is the only term
of that form.  By Lemma 2.3, if $W$ is good, then this g-poly
cannot be in a reduced form.  Hence, there is some non-trivial
relationship between the $p_{i}$, $q_{j}$. In other words, there
exists a subset, $P \subseteq \{p_{1},\ldots,p_{k}\}$, and a
subset, $Q \subseteq \{q_{1},\ldots,q_{k}\}$, such that $$p_{1} =
\sum\limits_{p \in P }{p } \ \ \ \text{and} \ \ \ q_{1} =
\sum\limits_{q \in Q }{q},$$ in which these equations do not both
represent the trivial relations $p_{1} = p_{1}$, $q_{1} = q_{1}$.
This proves the theorem. \qquad\qed
\enddemo

\

\subhead 4. Positive g-words\endsubhead We next give a strong
result for g-words with positive exponents, called {\it positive
g-words}. For a class $k$ g-word $W$ in standard form, consider
the list $L = L(W)$ of $2k$ real numbers: $$L = \{p_{1}+q_{1}, \
q_{1}+p_{2}, \ p_{2}+q_{2}, \ \ldots , \ p_{k}+q_{k}, \
q_{k}+p_{1}\},$$ the cyclically consecutive, pair-wise sums of the
exponents. Enumerate the elements of $L$ as $L_{1} = p_{1}+q_{1},\
L_{2} = q_{1}+p_{2}, \ \ldots , \ L_{2k} = q_{k}+p_{1}$.  Now,
suppose that the minimum of $L$ appears $m$ times, and let $L_{i_1
},L_{i_2 },\ldots,L_{i_m}$ be the appearances of this minimum in
$L$;  of course, $\{i_{1},\ldots,i_{m}\} \subseteq
\{1,\ldots,2k\}$.

Let $\#(L)_{\text{o}}$ denote the number of $i_{r}$ ($r =
1,\ldots,m$) that are odd, and let $\#(L)_{\text{e}}$ denote the
number of $i_{r}$ that are even. For example, the word $W =
A^{4}B^{2}A^{3}B^{3}A^{4}B^{1}$ has $L = \{6,5,6,7,5,5\}$, and so
$\#(L)_{\text{o}}$ = 1, $\#(L)_{\text{e}}$ = 2.  We call a word
{\it exact} if $\#(L)_{\text{o}} = \#(L)_{\text{e}}$. A word is
{\it inexact} if it is not exact.  We note combinatorially that
nearly symmetric g-words are exact, but that there are exact
g-words that are not nearly symmetric: $AB^{2}AB^{3}A^{4}B^{5}$.
In a sense, exact words are a first order (linear) approximation
to near symmetry.  Of course, a class 2 g-word is nearly symmetric
if and only if it is exact.  We then have the following.

\

\proclaim{Theorem 4.1}\rm All positive, good g-words are exact.
\endproclaim

\

We should remark at this point that Theorems 3.1 and 4.1 are very
different statements.  For example, consider the word $W =
A^{1}B^{2}A^{2}B^{2}A^{3}B^{4}$.  In this case, $$p_{1} = p_{1},\
q_{1} = q_{2}; \ \ p_{2} = p_{2},\ q_{2} = q_{1}; \ \text{and}\ \
p_{3} = p_{1}+p_{2},\ q_{3} = q_{3}$$ represent nontrivial pairs
of relations.  Thus, we cannot conclude that $W$ is bad using only
Theorem 3.1. However, from $L(W) = \{3,4,4,5,7,5\}$, it follows
that $\#(L)_{\text{o}} = 1, \#(L)_{\text{e}} = 0$, and so $W$ is
bad by Theorem 4.1.

\

Before we prove the theorem, we need a lemma that pertains only to
positive g-words. It allows us to find positive semidefinite $A$
and $B$ for which Tr$[W]$ is non-positive instead of finding PD
ones.  It will then be enough to show the theorem for the
following parameterization of the matrices $A,B$. As before let $S
= UU^{T}$ for some unitary matrix $U$, but this time assume $$S =
\left( {\matrix
   {s_{1,1} } & {s_{2,1} } & {s_{3,1} }  \\
   {s_{2,1} } & {s_{2,2} } & {s_{3,2} }  \\
   {s_{3,1} } & {s_{3,2} } & {s_{3,3} }  \\
 \endmatrix } \right), \ A = \left( {\matrix
   1 & 0 & 0  \\
   0 & x & 0  \\
   0 & 0 & 0  \\
 \endmatrix } \right), \ E = \left( {\matrix
   1 & 0 & 0  \\
   0 & y & 0  \\
   0 & 0 & 0  \\
 \endmatrix } \right), \ B = SE\overline S, \tag 4.1$$ for $x,y > 0$.

\

\proclaim{Lemma 4.2}\rm Suppose there exist $S$, $A$, and $E$ with
the parameterization (4.1) that give a positive g-word, $W$, a
nonzero imaginary part for its trace.  Then, $W$ is bad.
\endproclaim

\

\demo{Proof} Assume that we have found an $S$, $A$, and $E$ as
above that give the positive g-word, $W$, a nonzero imaginary part
for its trace.  Given $\varepsilon > 0$, let $W_{\varepsilon}$
denote the matrix product produced by replacing $A$ with
$A_{\varepsilon} =$ diag$(1,x,\varepsilon)$ and $E$ with
$E_{\varepsilon} =$ diag$(1,y,\varepsilon)$.  We can, therefore,
write Tr$[W_{\varepsilon}]$ as the trace of: $$ A_\varepsilon
^{p_1 } SE_\varepsilon ^{q_1 } \overline S A_\varepsilon  ^{p_2 }
SE_\varepsilon  ^{q_2 } \overline S  \cdots A_\varepsilon  ^{p_k }
SE_\varepsilon  ^{q_k } \overline S.$$  The imaginary part of this
product will be the same as that of $W$ except for an additional
(possibly 0) expression involving sums of positive powers of
$\varepsilon$.  Since the imaginary part of $W$ was assumed to be
nonzero, by continuity, we can choose $\varepsilon$ small enough
so that this imaginary part stays nonzero. \qquad\qed
\enddemo

\

The trace, Tr[$A^{p_1 } SE^{q_1 } \overline S \cdots A^{p_k }
SE^{q_k } \overline S$], of a word under the assumption of (4.1)
can be viewed as a g-poly in $x$, $y$, with exponents involving
the $p_{i}$ and $q_{j}$.  As before, it is not hard to see that
this trace will be a constant plus a (formal) sum of terms of the
form, $$c\,x^{p_{i_1 }  + p_{i_2 }  +  \cdots  + p_{i_u } }
y^{q_{j_1 }  + q_{j_2 }  + \cdots  + q_{j_v } },$$ in which $c \in
\Bbb C$, \{$p_{i_1 },\ldots ,p_{i_u } \} \subseteq \{p_{1},\ldots,
p_{k}\}$, and \{$q_{j_1 }, \ldots,q_{j_v } \} \subseteq
\{q_{1},\ldots, q_{k}\}$. Setting $y_{1} = y_{2} = 0$ in Lemma 3.3
above (which is a valid maneuver for positive g-words) gives us a
method to determine the coefficients of such terms, and an
argument similar to the one in Lemma 3.4 shows us that the
constant term in this trace is real.

We next compute the coefficients that are attached to certain
terms.  The aim is to discover which terms can be made to have
non-real coefficients.  This will give us insight into what values
of $p_{i}$ and $q_{j}$ guarantee that the word's trace can be made
non-real.  Our main result stems from the following.

\

\proclaim{Lemma 4.3}\rm Assuming (4.1), let $P = \{p_{i_1 },
\ldots ,p_{i_u }\} \subseteq \{p_{1},\ldots, p_{k}\}$ and let $Q =
\{q_{j_1 },\ldots ,q_{j_v }\} \subseteq \{q_{1},\ldots, q_{k}\}$.
Then, if the coefficient of the term $$ x^{p_{i_1 } + p_{i_2 }  +
\cdots + p_{i_u } } y^{q_{j_1 }  + q_{j_2 }  + \cdots  + q_{j_v }
} $$ is not real, there must be a $p_{i} \in P$ and a $q_{j} \in
Q$ such that $p_{i}$ and $q_{j}$ are adjacent.
\endproclaim

\

\demo{Proof} We will show the contrapositive.  Let $P = \{p_{i_1
},\ldots ,p_{i_u } \} \subseteq \{p_{1}, \ldots, p_{k}\}$ and let
$Q = \{ q_{j_1 },\ldots ,q_{j_v } \} \subseteq \{q_{1}, \ldots,
q_{k}\}$.  Then, we prove that if there are no adjacent
$p_{i},q_{j}$ in the $P,Q$, the coefficient of $$ x^{p_{i_1 }  +
p_{i_2 } + \cdots  + p_{i_u } } y^{q_{j_1 }  + q_{j_2 }  + \cdots
+ q_{j_v } }$$ is real.

If both $P$ and $Q$ are empty, then there is nothing to show (the
constant term is real).  Next, notice that it suffices to prove
the result under the assumption that $p_{1} \in P$.  This is
because if $P \neq \emptyset$, we can perform an appropriate
cycling (which will not change the coefficient) and relabel
variables. And if $P = \emptyset$, an interchange (which swaps
$x,y$) and a cycling will put our word into this form (all that
will change is the conjugacy of the coefficient by our
parameterization (4.1)). Using Lemma 3.3, it suffices to consider
the matrix $$J = A^{p_{i_1}} SF\overline S FSF\overline S \cdots
SE^{q_{j_1 } } \overline S FS \cdots SE^{q_{j_v } } \overline S
\cdots SF\overline S $$ where $F$ as in Lemma 3.3 replaces in $$
A^{p_1 } SE^{q_1 } \overline S A^{p_2 } SE^{q_2 } \overline S
\cdots A^{p_k } SE^{q_k } \overline S $$ each instance of $A^{p_i
}$ ($E^{q_j }$) in which $p_{i} \notin P$ ($q_{j} \notin Q$).  We
will perform an induction on the class number $k$ of a word,
utilizing a special form of the matrices $J$.

Assume that for some $k$ and for each $P,Q$ with no adjacent
elements, $J$ has the form $$ \left( {\matrix
   * & * & *  \\
   {s_{2,1} \overline {s_{1,1} } \,c_1  \cdot x^{p_{i_1 }
   + \cdots  + p_{i_u } } y^{q_{j_1 } +  \cdots  + q_{j_v } }  + g_1} & {c_2  \cdot x^{p_{i_1 }
    + \cdots  + p_{i_u } } y^{q_{j_1 } +  \cdots  + q_{j_v } }  + g_2} & *  \\
   0 & 0 & 0  \\
 \endmatrix } \right)
$$ where the constants $c_{1}$ and $c_{2}$ are real, $g_{1} =
g_{1}(x,y)$ and $g_{2} = g_{2}(x,y)$ are g-polys, and the term $$
x^{p_{i_1 } + \cdots + p_{i_u } } y^{q_{j_1 } +  \cdots + q_{j_v }
}$$ doesn't appear in the (1,1) or (1,2) locations in the above
matrix or in $g_{1}(x,y),g_{2}(x,y)$. Notice also that a matrix
with the form described above has a trace in which a desired term
has a real coefficient (namely, $c_{2}$).

The base case for our analysis will be $k = 1$.  In this
situation, the only term that doesn't violate the hypotheses about
adjacency is $x^{p_1 }$.  A calculation reveals that,

$$ A^{p_1 } SF\overline S = \left( {\matrix
   {s_{1,1} \overline {s_{1,1} } } & {s_{1,1} \overline {s_{2,1} } } & {s_{1,1} \overline {s_{3,1} } }  \\
   {x^{p_1 } s_{2,1} \overline {s_{1,1} } } & {x^{p_1 } s_{2,1} \overline {s_{2,1} } } & {x^{p_1 } s_{2,1} \overline {s_{3,1} } }  \\
   0 & 0 & 0  \\
 \endmatrix } \right).
$$

This matrix has the form given above and, thus, can be used as a
base case for our induction.

We now proceed with the induction.  Let $k$ be given in which the
lemma is true, and assume that we are given a class $k+1$ word and
a pair $P,Q$ with no adjacent elements.  Then, the coefficient of
the term, $$ x^{p_{i_1 } + \cdots  + p_{i_u } } y^{q_{j_1 } +
\cdots  + q_{j_v } },$$ can be found by examining the matrix
$$J_{1} = A^{p_{i_1 } } SF\overline S FSF\overline S  \cdots
SE^{q_{j_1 } } \overline S FS \cdots SE^{q_{j_{v - 1} } }
\overline S  \cdots SF\overline S $$ involving $k$ terms
(originally looking like $A^{p_1 } SE^{q_1 } \overline S \cdots
A^{p_k } SE^{q_k } \overline S$) multiplied by an appropriately
transformed (using Lemma 3.3) $A^{p_{k + 1} } SE^{q_{k + 1} }
\overline S$.

If $q_{k+1} \in Q$, then since $p_{1} \in P$, $p_{1}$ and
$q_{k+1}$ would be adjacent, contrary to our assumption.  Hence,
we may assume $q_{k+1} \notin Q$.  This leaves us with two cases:
$p_{k+1} \in P$ or $p_{k+1} \notin P$.

\

Case 1:  $p_{k+1} \notin P = \{ p_{i_1 },\ldots ,p_{i_u }\}$.

\

By the inductive hypothesis, the expression $J_{1}$ above has the
form: $$ \left( {\matrix
   * & * & *  \\
   {s_{2,1} \overline {s_{1,1} } \,c_1  \cdot x^{p_{i_1 }
   + \cdots  + p_{i_u } } y^{q_{j_1 } +  \cdots  + q_{j_v } }  + g_1} & {c_2  \cdot x^{p_{i_1 }
    + \cdots  + p_{i_u } } y^{q_{j_1 } +  \cdots  + q_{j_v } }  + g_2} & *  \\
   0 & 0 & 0  \\
 \endmatrix } \right).
$$ And since $p_{k+1} \notin P$, $q_{k+1} \notin Q$, we need only
concern ourselves with the product of $J_{1}$ with $FSF\overline
S$. From before, we know that $$ FSF\overline S \;\; = \;\;\left(
{\matrix
   {s_{1,1} \overline {s_{1,1} } } & {s_{1,1} \overline {s_{2,1} } } & {s_{1,1} \overline {s_{3,1} } }  \\
   0 & 0 & 0  \\
   0 & 0 & 0  \\
 \endmatrix } \right).
$$ Therefore, it is clear that $J_1 FSF\overline S$ preserves the
special form of the matrix discussed in the inductive hypothesis.

\

Case 2:  $p_{k+1} =  p_{i_u }  \in P$.

\

Since $p_{k+1} \in P$, we cannot have $q_{k} \in Q$ because then
$p_{k+1}$ and $q_{k}$ would be adjacent.  Therefore, we must have
that the product $$J_{1} = A^{p_{i_1 } } SF\overline S  \cdots
SE^{q_{j_1 } } \overline S FS \cdots SE^{q_{j_v } } \overline S
\cdots SF\overline S,$$ involving $k$ terms (originally looking
like $A^{p_1 } SE^{q_1 } \overline S \cdots A^{p_k } SE^{q_k }
\overline S$) complies with the induction hypothesis (the sets $P
\backslash \{p_{k+1}\}$ and $Q$ have no adjacent elements). Hence,
we need only look at $J_{1}$ multiplied by the matrix, $A^{p_{i_u
} } SF\overline S$, where by induction, $J_{1}$ looks like: $$
\left( {\matrix
   * & * & *  \\
   {s_{2,1} \overline {s_{1,1} } \,c_1  \cdot x^{p_{i_1 }
   + \cdots  + p_{i_{u-1} } } y^{q_{j_1 } +  \cdots  + q_{j_v } }  + g_1} & {c_2  \cdot x^{p_{i_1 }
    + \cdots  + p_{i_{u-1} } } y^{q_{j_1 } +  \cdots  + q_{j_v } }  + g_2} & *  \\
   0 & 0 & 0  \\
 \endmatrix } \right).
$$

A simple computation gives us

$$ A^{p_{i_u } } SF\overline S = \left( {\matrix
   {s_{1,1} \overline {s_{1,1} } } & {s_{1,1} \overline {s_{2,1} } } & {s_{1,1} \overline {s_{3,1} } }  \\
   {x^{p_{i_u } } s_{2,1} \overline {s_{1,1} } } & {x^{p_{i_u } } s_{2,1} \overline {s_{2,1} } } & {x^{p_{i_u } } s_{2,1} \overline {s_{3,1} } }  \\
   0 & 0 & 0  \\
 \endmatrix } \right).
$$

As before, it is not too difficult to see that the new matrix
produced, $ J_1 A^{p_{i_u } } SF\overline S $, will preserve the
special form needed for the induction.  This completes the
induction and the proof.   \qquad\qed
\enddemo

\

We have shown above a necessary condition on a term for it to have
a non-real coefficient.  We need one more lemma to guarantee that
there are some terms that that can be made to have non-real
coefficients.  This result is an analog of Lemma 3.5.

\

\proclaim{Lemma 4.4}\rm The coefficient of any term, $x^{p_{i_1 }
} y^{q_{j_1 } }$ in a class $k$ word ($k > 1$) in which $p_{i_1 }
$ and $q_{j_1 }$ are adjacent is given by: $$ \overline {s_{2,2} }
\left( {s_{2,1} } \right)^2 \overline {s_{1,1} } \left( {s_{1,1}
\overline {s_{1,1} } } \right)^{k - 2} \ \ \ \text{or} \ \ \
s_{2,2} \left( {\overline {s_{2,1} } } \right)^2 s_{1,1} \left(
{s_{1,1} \overline {s_{1,1} } } \right)^{k - 2}.$$
\endproclaim

\

\demo{Proof} We first notice that we can assume we are dealing
with the term $x^{p_1 } y^{q_1 }$ by cycling and (possibly) a
reversal.  Therefore, by Lemma 3.3, we need only compute $A^{p_1 }
SE^{q_1 } \overline S FSF\overline S  \cdots FSF\overline S$. This
is just $A^{p_1 } SE^{q_1 } \overline S \left( {FSF\overline S }
\right)^{k - 1}$.  A straightforward calculation gives us, $$
\left( {FSF\overline S } \right)^{k - 1}  = \left( {\matrix
   {\left( {s_{1,1} \overline {s_{1,1} } } \right)^{k - 1} } & {s_{1,1} \overline {s_{2,1} } \left( {s_{1,1} \overline {s_{1,1} } } \right)^{k - 2} } & {s_{1,1} \overline {s_{3,1} } \left( {s_{1,1} \overline {s_{1,1} } } \right)^{k - 2} }  \\
   0 & 0 & 0  \\
   0 & 0 & 0  \\
 \endmatrix } \right).
$$ Additionally, $A^{p_1 } SE^{q_1 } \overline S$ is $$\left(
{\matrix
   {s_{1,1} \overline {s_{1,1} }  + y^{q_1 } s_{2,1} \overline {s_{2,1} } } & {s_{1,1} \overline {s_{2,1} }  + y^{q_1 } s_{2,1} \overline {s_{2,2} } } & {* }  \\
   {x^{p_1} \left( {s_{2,1} \overline {s_{1,1} }  + y^{q_1 } s_{2,2} \overline {s_{2,1} } } \right)} & {x^{p_1 } \left( {s_{2,1} \overline {s_{2,1} }  + y^{q_1 } s_{2,2} \overline {s_{2,2} } } \right)} & {*}  \\
   0 & 0 & 0  \\
 \endmatrix } \right).$$ Therefore, taking the product of these two matrices and
computing the trace gives us a coefficient of $$s_{2,2} \left(
{\overline {s_{2,1} } } \right)^2 s_{1,1} \left( {s_{1,1}
\overline {s_{1,1} } } \right)^{k - 2} $$
 for the $x^{p_1 } y^{q_1 }$
 term as stated.   \qquad\qed
\enddemo

\

It should be clear that consecutive terms such as $x^{p_i } y^{q_i
}$ and $x^{p_{i + 1} } y^{q_i }$ have conjugate coefficients (by a
reversal and cycling) and that these coefficients can, in fact, be
made non-real (for example, using $U$ again from the proof of
Theorem 2.4).

\

We are now ready to prove Theorem 4.1.

\

\demo{Proof of Theorem 4.1} We will prove that an inexact word is
bad.  Take $A$, $S$, and $E$ as in (4.1).  As before, enumerate
the elements of $L$ as $L_{1} = p_{1}+q_{1},\ L_{2} = q_{1}+p_{2},
\ \ldots , \ L_{2k} = q_{k}+p_{1}$.  Now, suppose that the minimum
of $L$ appears $m$ times, and let $L_{i_1},L_{i_2},\ldots,L_{i_m}$
be the appearances of this minimum in $L$.  Without loss of
generality, we suppose it is the one involving the two terms
$p_{1}$ and $q_{1}$.

From Lemma 4.4, we know that the coefficients of these $m$ terms
can be made non-real for any class $k$ word's trace. From Lemma
4.3, however, we also know that any other term that appears in the
trace of $A^{p_1 } B^{q_1 } \cdots A^{p_k } B^{q_k }$ must contain
two adjacent $p_{i},q_{j}$ in the exponents of $x$ and $y$ in
order for it to have a non-real coefficient.

Set $x = y$, and notice that the imaginary part of the trace of an
inexact word cannot be zero for all $x > 0$. This is because there
are exactly $m$ terms of the same minimum degree, $p_{1}+q_{1}$
(since all of the $p_{i}$'s and $q_{j}$'s are positive), and they
all have nonzero coefficients. From the discussion after Lemma
4.4, it follows that the sum of these coefficients (whose signs
only alternate) is a non-zero constant times
$$\sum\limits_{k = 1}^m {\left( { - 1} \right)^{i_k }  = } \;\;\#
\left( L \right)_{\text{e}} - \# \left( L \right)_{\text{o}}.$$ If
the word is inexact, then by definition this sum will be non-zero.
Therefore, taking $x$ small enough, we can produce positive
semidefinite Hermitian $A$ and $B$ that give $W$ a non-real trace.
This gives us positive definite $A$ and $B$ by Lemma 4.2, and
concludes the proof of the theorem. \qquad\qed
\enddemo

\

If a positive g-word $W$ has $L(W)$ with $2k$ distinct elements,
then there is only one minimum element in $L$.  Hence, in this
case $W$ is inexact, and so Theorem 4.1 gives us the immediate
fact.

\

\proclaim{Corollary 4.5}\rm If $W$ is a positive g-word and if
$L(W)$ has distinct elements, then $W$ is bad.
\endproclaim

\

\subhead 5. Remark\endsubhead As a final remark, we note that we
can prove Conjecture 1.3 for words of class 3 and 4.  These proofs
contain arguments similar to those in Theorem 2.4, however, they
are much more cumbersome and do not shed any light on what is
happening in general.

\

\enddocument